\input amstex
\input amsppt.sty
\documentstyle{amsppt}
\pageheight{43pc}
\pagewidth{28pc}
\NoBlackBoxes
\TagsOnRight \pageno=1 \nologo
\def\Z{\Bbb Z}
\def\N{\Bbb N}

\def\l{\left}
\def\r{\right}
\def\bg{\bigg}
\def\({\bg(}
\def\[{\bg\lfloor}
\def\){\bg)}
\def\]{\bg\rfloor}
\def\t{\text}
\def\f{\frac}

\def\bi{\binom}
\def\eq{\equiv}

\def\ls{\leqslant}

\def\mo{\roman{mod}}

\def\Proof{\noindent{\it Proof}}

\def\Remark{\medskip\noindent{\it  Remark}}

\def\Ack{\medskip\noindent {\bf Acknowledgment}}
\hbox {Illinois J. Math. 56(2012), no.\,3, 967--979.}
\bigskip
\topmatter
\title A refinement of a congruence result by van Hamme and Mortenson\endtitle
\author Zhi-Wei Sun\endauthor
\leftheadtext{Zhi-Wei Sun} \rightheadtext{A refinement of a
congruence result}
\affil Department of Mathematics, Nanjing University\\
 Nanjing 210093, People's Republic of China
  \\  zwsun\@nju.edu.cn
  \\ {\tt http://math.nju.edu.cn/$\sim$zwsun}
\endaffil
\abstract Let $p$ be an odd prime. In 2008 E. Mortenson proved van Hamme's following conjecture:
$$\sum_{k=0}^{(p-1)/2}(4k+1)\bi{-1/2}k^3\eq (-1)^{(p-1)/2}p\ \ (\mo\ p^3).$$
In this paper we show further that
$$\align\sum_{k=0}^{p-1}(4k+1)\bi{-1/2}k^3\eq&\sum_{k=0}^{(p-1)/2}(4k+1)\bi{-1/2}k^3
\\\eq& (-1)^{(p-1)/2}p+p^3E_{p-3}
\ (\mo\ p^4),\endalign$$
where $E_0,E_1,E_2,\ldots$ are Euler numbers. We also prove that if $p>3$ then
$$\sum_{k=0}^{(p-1)/2}\f{20k+3}{(-2^{10})^k}\bi{4k}{k,k,k,k}\eq(-1)^{(p-1)/2}p(2^{p-1}+2-(2^{p-1}-1)^2)\ (\mo\ p^4).$$
\endabstract
\thanks 2010 {\it Mathematics Subject Classification}.\,Primary 11B65;
Secondary 05A10,  11A07, 11B68.
\newline\indent Supported by the National Natural Science
Foundation (grant 11171140) of China and and the Priority Academic
Program Development of Jiangsu Higher Education Institutions.
\endthanks
\endtopmatter
\document

\heading{1. Introduction}\endheading

In 1859 G. Bauer obtained the identity
$$\sum_{k=0}^\infty(4k+1)\bi{-1/2}k^3=\f2{\pi}$$
which was later reproved by S. Ramanujan [R] in 1914. (Note that
$\bi{-1/2}k=\bi{2k}k/(-4)^k$ for all $k=0,1,2,\ldots$.) In 1997 van
Hamme [vH] conjectured that
$$\sum_{k=0}^{p-1}(4k+1)\bi{-1/2}k^3=\sum_{k=0}^{p-1}(4k+1)\f{\bi{2k}k^3}{(-64)^k}\eq (-1)^{(p-1)/2}p\ (\mo\ p^3)
$$ for any odd prime $p$,
which was first confirmed by E. Mortenson [Mo] in 2008 via a deep
method involving the $p$-adic $\Gamma$-function and Gauss and Jacobi
sums.

Throughout this paper, for an odd prime $p$, we use $(\f\cdot{p})$
to denote the Legendre symbol. Recall that the Euler numbers
$E_0,E_1,E_2,\ldots$ are integers given by
$$E_0=1\ \ \t{and}\ \ \sum^n\Sb k=0\\2\mid k\endSb \bi nk E_{n-k}=0\ \ (n=1,2,3,\ldots).$$
It is well known that
$$\f{2e^x}{e^{2x}+1}=\sum_{n=0}^\infty E_n\f{x^n}{n!}\qquad \t{for}\ |x|<\f{\pi}2.$$

In this paper we obtain the following refinement of the congruence
by van Hamme and Mortenson via an elementary approach.

\proclaim{Theorem 1.1} Let $p$ be an odd prime. Then
$$\sum_{k=0}^{p-1}(4k+1)\f{\bi{2k}k^3}{(-64)^k}
\eq\sum_{k=0}^{(p-1)/2}(4k+1)\f{\bi{2k}k^3}{(-64)^k}\eq p\l(\f{-1}p\r)+p^3E_{p-3}\
(\mo\ p^4).\tag1.1$$
\endproclaim
\Remark\ 1.1. The only previously proved congruence mod $p^4$ of the same kind is the following one
conjectured by van Hamme [vH] and confirmed by L. Long [Lo]:
$$\sum_{k=0}^{(p-1)/2}(6k+1)\f{\bi{2k}k^3}{256^k}\eq p\l(\f{-1}p\r)\ \ (\mo\ p^4)\quad\t{for any prime}\ p>3.$$

For each nonnegative integer $k$, it is clear that
$$\bi{4k}{k,k,k,k}=\f{(4k)!}{k!^4}=\bi{4k}{2k}\bi{2k}k^2.$$
In a way similar to the proof of Theorem 1.1, we also deduce the following result.

 \proclaim{Theorem 1.2} Let $p>3$ be a prime. Then
$$\sum_{k=0}^{(p-1)/2}\f{20k+3}{(-2^{10})^k}\bi{4k}{k,k,k,k}\eq p\l(\f{-1}p\r)\l(2^{p-1}+2-(2^{p-1}-1)^2\r)\
(\mo\ p^4).\tag1.2$$
\endproclaim
\Remark.\ 1.2. (a) The congruence in Theorem 1.2 gives the mod $p^4$ analogy of the Ramanujan series
$$\sum_{k=0}^\infty\f{20k+3}{(-2^{10})^k}\bi{4k}{k,k,k,k}=\f 8{\pi}.$$
See [BB], [BBC] and [Be,\ pp.353-354] for more such series. The mod
$p^3$ analogy of the above series is known (cf. [Zu]).

(b) By the same method, the author ever proved that
$$\sum_{k=0}^{p-1}\f{20k+3}{(-2^{10})^k}\bi{4k}{k,k,k,k}\eq3p\l(\f{-1}p\r)+3p^3E_{p-3}\ (\mo\ p^4)\tag1.3$$
for any odd prime $p$; unfortunately he has lost the draft containing the complicated details.

\medskip

Theorems 1.1 and 1.2 will be proved in Sections 2 and 3 respectively.

The author [Su2, Conjecture 5.1] raised several conjectures similar to (1.1).
Here we pose a new conjecture motivated by the Ramanujan series
$$\sum_{k=0}^\infty\f{7k+1}{648^k}\bi{4k}{k,k,k,k}=\f 9{2\pi}.$$

\proclaim{Conjecture 1.1} For any prime $p>3$ we have
$$\sum_{k=0}^{p-1}\f{7k+1}{648^k}\bi{4k}{k,k,k,k}\eq p\l(\f{-1}p\r)-\f 53p^3E_{p-3}\ (\mo\ p^4).\tag1.4$$
Also, for $n=2,3,\ldots$ we have
$$\f1{2n(2n+1)\bi{2n}n}\sum_{k=0}^{n-1}(7k+1)\bi{4k}{k,k,k,k}648^{n-1-k}\in\Z$$
unless $2n+1$ is a power of $3$ in which case the quotient is a rational number with denominator $3$.
\endproclaim
\Remark\ 1.3. It seems that the method for our proofs of (1.1) and
(1.2) does not work for (1.4).
\medskip

In 2010, L. L. Zhao, H. Pan and the author [ZPS] proved that
$$\sum_{k=1}^{p-1}\f{2^k}k\bi{3k}k\eq0\ \ (\mo\ p)$$ for any odd prime
$p$. Here we raise a further conjecture.

\proclaim{Conjecture 1.2} Let $p$ be an odd prime. Then
$$\sum_{k=1}^{p-1}\f{2^k}k\bi{3k}k\eq-\f3p(2^{p-1}-1)^2\ \ (\mo\ p^2)\tag1.5$$
and
$$\sum_{k=1}^{p-1}\f{2^k}{k^2}\bi{3k}k\eq6\l(\f{-1}p\r)E_{p-3}\ \ (\mo\ p).\tag1.6$$
Also,
$$p\sum_{k=1}^{p-1}\f1{k2^k\bi{3k}k}\eq\cases0\ (\mo\ p^2)&\t{if}\ p\eq1\ (\mo\ 4),
\\-3/5\ (\mo\ p^2)&\t{if}\ p\eq3\ (\mo\ 4),\endcases\tag1.7$$
$$p\sum_{k=1}^{p-1}\f1{k^22^k\bi{3k}k}\eq\f{1-4^{p-1}}{4p}\ (\mo\ p^2)\ \ \ \t{if}\ p>3,\tag1.8$$
and
$$\sum_{k=1}^{p-1}2^k\bi{3k}k\sum_{j=1}^k\f1{j^2}\eq0\ (\mo\ p)\ \ \ \t{if}\ p>5\ \t{and}\ p\eq1\ (\mo\ 4).\tag1.9$$
\endproclaim

\heading{2. Proof of Theorem 1.1}\endheading

We need some classical congruences.

\proclaim{Lemma 2.1} Let $p>3$ be a prime.

{\rm (i) (J. Wolstenholme [W])} We have
$$\sum_{k=1}^{p-1}\f1k\eq0\ (\mo\ p^2),\ \ \sum_{k=1}^{p-1}\f1{k^2}\eq0\ (\mo\ p),\tag2.1$$
and
$$\f12\bi{2p}p=\bi{2p-1}{p-1}\eq 1\ (\mo\ p^3).\tag2.2$$

{\rm (ii) (F. Morley [M])} We have
$$\bi{p-1}{(p-1)/2}\eq(-1)^{(p-1)/2} 4^{p-1}\ (\mo\ p^3).\tag2.3$$
\endproclaim

The most crucial lemma we need is the following sophisticated result.

\proclaim{Lemma 2.2 {\rm (Sun [Su1])}} Let $p$ be an odd prime. Then
$$\sum_{k=1}^{(p-1)/2}\f{4^k}{(2k-1)\bi{2k}k}\eq E_{p-3}-1+\l(\f{-1}p\r)\ (\mo\ p)\tag2.4$$ and
$$\sum_{k=1}^{(p-1)/2}\f{4^k}{k(2k-1)\bi{2k}k}\eq 2E_{p-3}\ (\mo\ p).\tag2.5$$
\endproclaim
\Remark\ 2.1. Actually  (2.4) and (2.5) are equivalent since
$$\f12\sum_{k=1}^n\f{4^k}{k\bi{2k}k}=\f{4^n}{\bi{2n}n}-1;$$
they are (1.3) and (3.1) of Sun [Su1] respectively.

\medskip
\noindent{\it Proof of Theorem 1.1}. (i) Clearly the first congruence in (1.1) has the following equivalent form:
$$\sum_{p/2<k<p}(4k+1)\f{\bi{2k}k^3}{(-64)^k}\eq0\ (\mo\ p^4).$$
For $k\in\{1,\ldots,(p-1)/2\}$, it is obvious that
$$\align\f1p\bi{2(p-k)}{p-k}=&\f1p\times\f{p!\prod_{s=1}^{p-2k}(p+s)}{((p-1)!/\prod_{0<t<k}(p-t))^2}
\\\eq&\f{(k-1)!^2}{(p-1)!/(p-2k)!}\eq-\f{(k-1)!^2}{(2k-1)!}=-\f2{k\bi{2k}k}\ (\mo\ p).
\endalign$$
(See also [Su2, Lemma 2.1].)
Thus
$$\align&\f1{p^3}\sum_{p/2<k<p}(4k+1)\f{\bi{2k}k^3}{(-64)^k}
\\=&\sum_{k=1}^{(p-1)/2}\f{4(p-k)+1}{(-64)^{p-k}}\l(\f{\bi{2(p-k)}{p-k}}p\r)^3
\\\eq&\sum_{k=1}^{(p-1)/2}(1-4k)(-64)^{k-1}\l(\f{-2}{k\bi{2k}k}\r)^3
\\=&-\f18\sum_{k=1}^{(p-1)/2}\f{4k-1}{k^3\bi{-1/2}k^3}=\sum_{k=1}^{(p-1)/2}\f{4k-1}{\bi{-3/2}{k-1}^3}
\\\eq&\sum_{k=0}^{(p-3)/2}\f{4(k+1)-1}{\bi{(p-3)/2}{k}^3}
=\f12\sum_{k=0}^{(p-3)/2}\f{(4k+3)+4((p-3)/2-k)+3}{\bi{(p-3)/2}{k}^3}
\\\eq&0\ (\mo\ p)
\endalign$$
and hence the first congruence in (1.1) follows.
\medskip

(ii) Below we prove the second congruence in (1.1). For $k,n=0,1,2,\ldots$ define
$$F(n,k)=\f{(-1)^{n+k}(4n+1)}{4^{3n-k}}\bi{2n}n^2\f{\bi{2n+2k}{n+k}\bi{n+k}{2k}}{\bi{2k}k}$$
and
$$G(n,k)=\f{(-1)^{n+k}(2n-1)^2\bi{2n-2}{n-1}^2}{2(n-k)4^{3(n-1)-k}}\bi{2(n-1+k)}{n-1+k}\f{\bi{n-1+k}{2k}}{\bi{2k}k}.$$
Clearly $F(n,k)=G(n,k)=0$ if $n<k$. It can be easily verified that
$$F(n,k-1)-F(n,k)=G(n+1,k)-G(n,k)$$
for all nonnegative integers $n$ and $k>0$ as observed by S. B.
Ekhad and D. Zeilberger [EZ].

Let $m=(p-1)/2$. In the spirit of the WZ (Wilf-Zeilberger) method (see the book of M. Petkov\v sek, H.
S. Wilf and D. Zeilberger [PWZ], and  [AZ] and [Z] for this
method), we have
$$\align \sum_{n=0}^mF(n,0)-F(m,m)=&\sum_{n=0}^mF(n,0)-\sum_{n=0}^mF(n,m)
\\=&\sum_{k=1}^m\(\sum_{n=0}^mF(n,k-1)-\sum_{n=0}^mF(n,k)\)
\\=&\sum_{k=1}^m\sum_{n=0}^m(G(n+1,k)-G(n,k))
=\sum_{k=1}^mG(m+1,k),
\endalign$$
that is,
$$\aligned&\sum_{n=1}^m\f{4n+1}{(-64)^n}\bi{2n}n^3-\f{4m+1}{4^{2m}}\bi{4m}{2m}\bi{2m}m
\\=&\sum_{k=1}^m\f{(-1)^{m+k+1}(2m+1)^2\bi{2m}m^2}{2(m+1-k)4^{3m-k}}\bi{2m+2k}{m+k}\f{\bi{m+k}{2k}}{\bi{2k}k}.
\endaligned\tag2.6$$

For $0<k\ls m=(p-1)/2$, clearly
$$\align&\f1p\bi{2m+2k}{m+k}=\f{(p-1)!(p+1)\cdots(p+2k-1)}{m!^2\prod_{j=1}^k((p+2j-1)/2)^2}
\\\eq&(-1)^{(p-1)/2}\f{(p-1)!}{\prod_{k=1}^{(p-1)/2}k(p-k)}\cdot\f{(2k-1)!}{((2k-1)!!/2^k)^2}
\\\eq&\l(\f{-1}p\r)\f{(2k-1)!}{((2k)!/(k!4^k))^2}=\l(\f{-1}p\r)\f{4^{2k}}{2k\bi{2k}k}\ (\mo\ p)
\endalign$$
and
$$\align\bi{m+k}{2k}\eq&\bi{k-1/2}{2k}=\f{\prod_{j=1}^k(-(2j-1)/2)(2j-1)/2}{(2k)!}
\\=&\f{(-1)^k((2k-1)!!)^2}{4^k(2k)!}=\f{((2k)!/\prod_{j=1}^k(2j))^2}{(-4)^k(2k)!}=\f{\bi{2k}k}{(-16)^k}\ (\mo\ p).
\endalign$$
Note also that
$$(4m+1)\bi{4m}{2m}=(2p-1)\bi{2p-2}{p-1}=p\bi{2p-1}p\eq p\ (\mo\ p^4)$$
by the Wolstenholme congruence (2.2). Thus, in view of the above and Morley's congruence (2.3), we obtain from
(2.6) that
$$\align&\sum_{k=0}^m(4k+1)\f{\bi{2k}k^3}{(-64)^k}-p(-1)^{(p-1)/2}
\\\eq& p^3\sum_{k=1}^m\f{(-1)^{k-1}4^{2k}}{2((p+1)/2-k)2^{3(p-1)-2k}2k\bi{2k}k(-16)^k}
\\\eq&\f{p^3}2\sum_{k=1}^{(p-1)/2}\f{4^k}{k(2k-1)\bi{2k}k}\ (\mo\ p^4)
\endalign$$
Combining this with (2.5) we get the second congruence in (1.1).

The proof of Theorem 1.1 is now complete. \qed

\heading{3. Proof of Theorem 1.2}\endheading

\proclaim{Lemma 3.1} Let $p$ be an odd prime. Then
$$\bi{(p-1)/2+k}{2k}\eq\f{\bi{2k}k}{(-16)^k}\ (\mo\ p^2).\tag3.1$$
\endproclaim
\Remark\ 3.1. (3.1) is easy, see [S, Lemma 2.2] for a proof.
\medskip

Recall that the harmonic numbers are those rational numbers
$$H_n:=\sum_{k=1}^n\f1k\ \ (n=1,2,\ldots),$$
together with $H_0=0$. For an odd prime $p$ we write $q_p(2)$ for
the Fermat quotient $(2^{p-1}-1)/p$.

\proclaim{Lemma 3.2 {\rm (E. Lehmer [L])}} For any odd prime $p$ we have
$$H_{(p-1)/2}\eq-2q_p(2)+p\,q_p(2)^2\ (\mo\ p^2).\tag3.2$$
\endproclaim

\proclaim{Lemma 3.3} Let $p$ be an odd prime. Then
$$\sum_{k=1}^{(p-1)/2}\f{H_{k-1}}k\eq2q_p(2)^2\ (\mo\ p).\tag3.3$$
\endproclaim
\Proof. For $k=1,\ldots,p-1$ we have
$$\f{\bi pk}p=\f{\bi{p-1}{k-1}}k=\f{(-1)^{k-1}}k\prod_{0<j<k}\l(1-\f pj\r)\eq\f{(-1)^{k-1}}k(1-pH_{k-1})\ (\mo\ p^2).$$
Thus
$$\sum_{k=1}^{(p-1)/2}\f{pH_{k-1}-1}k\eq\f1p\sum_{k=1}^{(p-1)/2}(-1)^k\bi pk\pmod{p^2}.$$
As $\sum_{k=0}^{(p-1)/2}(-1)^k\bi pk$ is the coefficient of $x^{(p-1)/2}$ in  $(1-x)^p(1-x)^{-1}$,
we have
$$\align \f1p\sum_{k=1}^{(p-1)/2}(-1)^k\bi pk=\f{\bi{p-1}{(p-1)/2}(-1)^{(p-1)/2}-1}p
\eq\f{4^{p-1}-1}p\ (\mo\ p^2)
\endalign$$
with the help of Morley's congruence (2.3). Therefore, in view of Lehmer's congruence (3.2), we have
$$\align p\sum_{k=1}^{(p-1)/2}\f{H_{k-1}}k\eq&H_{(p-1)/2}+\f{2^{p-1}-1}p(2^{p-1}+1)
\\\eq&-2q_p(2)+p\,q_p(2)^2+q_p(2)(2+p\,q_p(2))
\\=&2p\,q_p(2)^2\ \ (\mo\ p^2)
\endalign$$
and hence (3.3) holds. \qed

\proclaim{Lemma 3.4} Let $p=2m+1$ be an odd prime. Then
$$\f{6m+1}{2^{8m}}\bi{6m}{3m}\bi{3m}m\eq p\l(\f{-1}p\r)\ (\mo\ p^4).\tag3.4$$
\endproclaim
\Proof. Observe that
$$\align &(6m+1)\bi{6m}{3m}\bi{3m}m=\f{(3m+1)\cdots(6m+1)}{m!(2m)!}
\\=&\f{(p+(p-1)/2)\cdots2p\cdots(3p-2)}{(p-1)!((p-1)/2)!}
=\f{(p+(p+1)/2)\cdots2p\cdots(3p-1)}{2\times(p-1)!((p-1)/2)!}
\\=&p\prod_{k=1}^{(p-1)/2}\f{(2p-k)(2p+k)}{k^2}\times\prod_{p/2<j<p}\f{2p+j}j
\\=&p(-1)^{(p-1)/2}\prod_{k=1}^{(p-1)/2}\l(1-\f{4p^2}{k^2}\r)\prod_{p/2<j<p}\l(1+\f{2p}j\r).
\endalign$$
Clearly
$$\prod_{k=1}^{(p-1)/2}\l(1-\f{4p^2}{k^2}\r)\eq1-4p^2\sum_{k=1}^{(p-1)/2}\f1{k^2}\eq1\ (\mo\ p^3)$$
since
$$2\sum_{k=1}^{(p-1)/2}\f1{k^2}\eq\sum_{k=1}^{(p-1)/2}\l(\f1{k^2}+\f1{(p-k)^2}\r)\eq\sum_{k=1}^{p-1}\f1{k^2}\eq0\ (\mo\ p).$$
So it suffices to prove that
$$\prod_{p/2<j<p}\l(1+\f{2p}j\r)\eq 2^{4(p-1)}\ (\mo\ p^3).\tag3.5$$

Observe that
$$\align&\prod_{p/2<j<p}\l(1+\f{2p}j\r)
\\\eq&1+2p\sum_{p/2<j<p}\f1j+4p^2\sum_{p/2<i<j<p}\f1{ij}
\\\eq&1+2p(H_{p-1}-H_{(p-1)/2})+2p^2\(\(\sum_{p/2<k<p}\f1k\)^2-\sum_{p/2<k<p}\f1{k^2}\)
\\\eq&1-2pH_{(p-1)/2}+2p^2(-H_{(p-1)/2})^2\qquad (\t{by}\ (2.1))
\\\eq&1-2p(p\,q_p(2)^2-2q_p(2))+2p^24q_p(2)^2\qquad(\t{by}\ (3.2))
\\=&1+4p\,q_p(2)+6p^2q_p(2)^2\eq(1+p\,q_p(2))^4=2^{4(p-1)}\ (\mo\ p^3).
\endalign$$
This proves (3.5) and hence (3.4) follows. \qed

\medskip
\noindent{\it Proof of Theorem 1.2}. (i) For $n,k\in\N$, define
$$F(n,k):=\f{(-1)^{n+k}(20n-2k+3)}{4^{5n-k}}\cdot\f{\bi{2n}n\bi{4n+2k}{2n+k}\bi{2n+k}{2k}\bi{2n-k}n}{\bi{2k}k}.
$$
and
$$G(n,k)
:=\f{(-1)^{n+k}}{4^{5n-4-k}}\cdot\f{n\bi{2n-1}{n-1}\bi{4n+2k-2}{2n+k-1}\bi{2n+k-1}{2k}\bi{2n-k-1}{n-1}}{\bi{2k}k}.$$
Clearly $F(n,k)=0$ if $n<k$. It can be easily verified that
$$F(n,k-1)-F(n,k)=G(n+1,k)-G(n,k)$$
for all nonnegative integers $n$ and $k>0$; the WZ-pair $F$ and $G$
stated in [Zu] was found in the spirit of [EZ] and [PWZ].

As in the proof of Theorem 1.1, for any positive integer $N$ we have
$$\sum_{n=0}^{N}F(N,0)-F(N,N)=\sum_{k=1}^{N}G(N+1,k),$$
that is,
$$\aligned&\sum_{n=0}^{N}\f{20n+3}{(-2^{10})^n}\bi{2n}n^2\bi{4n}{2n}
-\f{18N+3}{2^{8N}}\bi{6N}{3N}\bi{3N}{N}
\\=&(N+1)\bi{2N+1}{N}\sum_{k=1}^{N}\f{(-1)^{N+k+1}\bi{4N+2k+2}{2N+k+1}\bi{2N+k+1}{2k}\bi{2N-k+1}{N}}{4^{5(N+1)-4-k}\bi{2k}k}.
\endaligned\tag3.6$$
For $1\ls k\ls N$, clearly
$$\align&\bi{4N+2k+2}{2N+k+1}\bi{2N+k+1}{2k}\bi{2N-k+1}N
\\=&\bi{4N+2k+2}{2k}\bi{4N+2}{2N-k+1}\bi{2N-k+1}N
\\=&\bi{4N+2k+2}{2k}\bi{4N+2}N\bi{3N+2}{N-k+1}.
\endalign$$
So we also have
$$\aligned&\sum_{n=0}^{N}\f{20n+3}{(-2^{10})^n}\bi{2n}n^2\bi{4n}{2n}
-\f{18N+3}{2^{8N}}\bi{6N}{3N}\bi{3N}{N}
\\=&(N+1)\bi{2N+1}{N}\bi{4N+2}N
\sum_{k=1}^{N}\f{(-1)^{N+k+1}\bi{4N+2k+2}{2k}\bi{3N+2}{N-k+1}}{4^{5N+1-k}\bi{2k}k}.
\endaligned\tag3.7$$

(ii) Let $m=(p-1)/2$. Observe that
$$(m+1)\bi{2m+1}{m}=p\bi{p-1}{(p-1)/2}\eq p(-1)^m4^{p-1}\ (\mo\ p^4)$$
by Morley's congruence (2.3). Also,
$$\align\bi{4m+2}m=&\bi{2p}{(p-1)/2}=\f{4p}{p+1}\bi{2p-1}p\bi{p-1}{(p-1)/2}\prod_{k=1}^{(p+1)/2}\l(1+\f pk\r)^{-1}
\\\eq&\f{4p}{p+1}(-1)^{(p-1)/2}4^{p-1}\prod_{k=1}^{(p+1)/2}\l(1-\f pk\r)
\\\eq&p4^p(-1)^m(1-p)(1-pH_{(p+1)/2})
\\\eq& p4^p (-1)^m(1-p)(1-2p+2p\,q_p(2))
\\\eq&p4^p(-1)^m(1-3p+2p\,q_p(2))\ (\mo\ p^3)
\endalign$$
by Lehmer's congruence (3.2). Therefore
$$\align \f{(m+1)\bi{2m+1}{m}\bi{4m+2}m}{4^{5m+1}}\eq&p^2\f{4^{2(p-1)}(1-3p+2p\,q_p(2))}{4^{4m}(1+p\,q_p(2))}
\\\eq&p^2(1-p\,q_p(2))(1-3p+2p\,q_p(2))
\\\eq&p^2(1-3p+p\,q_p(2))\ (\mo\ p^4).
\endalign$$
Observe that
$$\align&\sum_{k=1}^m(-1)^k\f{\bi{4m+2k+2}{2k}\bi{3m+2}{m-k+1}}{4^{-k}\bi{2k}k}
\\\eq&\sum_{k=1}^m(-1)^k\f{\bi{2p+2k}{2k}\bi{p+(p+1)/2}{(p+1)/2-k}}{4^{-k}\bi{(p-1)/2+k}{2k}(-16)^k}
\\=&\sum_{k=1}^m\f{(2p+1)\cdots(2p+2k)(p+k+1)\cdots(p+(p+1)/2)}{((p+1)/2-k)!4^k((p-1)/2+k)!/((p-1)/2-k)!}
\\=&\f{(p+1)\cdots(p+(p+1)/2)}{((p-1)/2)!}\sum_{k=1}^m\f{\prod_{j=1}^k(2p+2j-1)}{((p+1)/2-k)2^k\prod_{j=1}^k((p-1)/2+j)}
\\=&\f{3p+1}2\prod_{j=1}^{(p-1)/2}\l(1+\f pj\r)\sum_{k=1}^m\f{\prod_{j=1}^k(1+p/(p+2j-1))}{(p+1)/2-k}\ (\mo\ p^2)
\endalign$$
and hence
$$\align&\sum_{k=1}^m(-1)^k\f{\bi{4m+2k+2}{2k}\bi{3m+2}{m-k+1}}{4^{-k}\bi{2k}k}
\\\eq&\f{3p+1}2(1+pH_{(p-1)/2})\sum_{s=1}^m\f{1+p\sum_{j=1}^{(p+1)/2-s}1/(2j-1)}s
\\\eq&\f{1+3p-2p\,q_p(2)}2\(H_m+\sum_{s=1}^m\f ps\sum_{t=s}^{(p-1)/2}\f1{2((p+1)/2-t)-1}\)
\\\eq&\f{1+3p-2p\,q_p(2)}2\(H_m-\f p2\sum_{s=1}^m\f {H_m-H_{s-1}}s\)
\\\eq&\f{1+3p-2p\,q_p(2)}2\(H_m-\f p2H_m^2+\f p2\sum_{k=1}^m\f{H_{k-1}}k\)\ (\mo\ p^2).
\endalign$$
Applying Lemmas 3.2 and 3.3 we get
$$\align&\sum_{k=1}^m(-1)^k\f{\bi{4m+2k+2}{2k}\bi{3m+2}{m-k+1}}{4^{-k}\bi{2k}k}
\\\eq&\f{1+3p-2p\,q_p(2)}2\l(-2q_p(2)+p\,q_p(2)^2-\f p2\cdot4q_p(2)^2+\f p2\cdot 2q_p(2)^2\r)
\\\eq&-q_p(2)(1+3p-2p\,q_p(2))\ (\mo\ p^2).
\endalign$$

Let $L$ and $R$ denote the left-hand side and the right-hand side of (3.7) with $N=m$ respectively. By the above,
$$\align R\eq& p^2(1-3p+p\,q_p(2))(-1)^{m+1}(-q_p(2))(1+3p-2pq_p(2))
\\\eq& p^2(-1)^mq_p(2)(1-p\,q_p(2))
\\=&p\l(\f{-1}p\r)(2^{p-1}-1)(1-(2^{p-1}-1))\ (\mo\ p^4).
\endalign$$
On the other hand, with the help of Lemma 3.4 we have
$$L=\sum_{k=0}^{(p-1)/2}\f{20k+3}{(-2^{10})^k}\bi{4k}{k,k,k,k}-3p\l(\f{-1}p\r)\ (\mo\ p^4).$$
So (3.7) with $N=m$ yields the desired (1.2). We are done. \qed

\Ack. The author is grateful to the referee for helpful comments.

\enddocument


\widestnumber\key{PWZ}

 \Refs

 \ref\key AZ\by T. Amdeberhan and D. Zeilberger \paper Hypergeometric series acceleration via the WZ method
\jour Electron. J. Combin.\vol4\yr 1997\pages no.\,2, \#R3\endref

\ref\key BB\by N. D. Baruah and B. C. Berndt\paper Eisenstein series
and Ramanujan-type series for $1/\pi$ \jour Ramanujan J.\vol 23\yr
2010\pages 17--44\endref

\ref\key BBC \by N. D. Baruah, B. C. Berndt and H. H. Chan\paper
Ramanujan's series for $1/\pi$: a survey\jour Amer. Math.
Monthly\vol 116\yr 2009\pages 567--587\endref

\ref\key Be\by B. C. Berndt\book Ramanujan's Notebooks, Part IV\publ
Springer, New York, 1994\endref

\ref\key EZ\by S. B. Ekhad and D. Zeilberger\paper A WZ proof of
Ramanujan's formula for $\pi$\jour in: Geometry, Analysis, and
Mechanics (J. M. Rassias, ed.), World Sci. Publ., Singapore, 1994,
107--108\endref

\ref\key L\by E. Lehmer\paper On congruences involving Bernoulli
numbers and the quotients of Fermat and Wilson\jour Ann. of Math.
(2)\vol 39\yr 1938\pages 350--360\endref

\ref\key Lo\by L. Long \paper Hypergeometric evaluation identities
and supercongruences \jour Pacific J. Math.\vol 249\yr 2011\pages
405--418\endref

\ref\key M\by F. Morley\paper Note on the congruence
$2^{4n}\equiv(-1)^n(2n)!/(n!)^2$, where $2n+1$ is a prime\jour Ann.
of Math. \vol 9\yr 1895\pages 168--170\endref

\ref\key Mo\by E. Mortenson\paper A $p$-adic supercongruence conjecture of van Hamme
\jour Proc. Amer. Math. Soc.\vol 136\yr 2008\pages 4321--4328\endref

\ref\key PWZ\by M. Petkov\v sek, H. S. Wilf and D. Zeilberger\book $A=B$ \publ A K Peters, Wellesley, 1996\endref

\ref\key R\by S. Ramanujan\paper Modular equations and approximations to $\pi$
\jour Quart. J. Math. (Oxford) (2)\vol45\yr1914
\pages 350--372\endref

\ref\key S\by Z. H. Sun\paper Congruences concerning Legendre
polynomials \jour Proc. Amer. Math. Soc. \vol 139\yr 2011\pages 1915--1929\endref

\ref\key Su1\by Z. W. Sun\paper On congruences related to central
binomial coefficients \jour J. Number Theory \vol 131\yr
2011\pages 2219-2238\endref

\ref\key Su2\by Z. W. Sun\paper Supper congruences and Euler numbers
\jour Sci. China Math. \vol 54\yr 2011\pages 2509--2535\endref

\ref\key vH\by L. van Hamme\paper Some conjectures concerning partial sums of generalized hypergeometric series
\jour in: $p$-adic Functional Analysis (Nijmegen, 1996), pp. 223--236,
Lecture Notes in Pure and Appl. Math., Vol., 192, Dekker, 1997\endref

\ref\key W\by J. Wolstenholme\paper On certain properties of prime numbers\jour Quart. J. Appl. Math.
\vol 5\yr 1862\pages 35--39\endref

\ref\key Z\by D. Zeilberger\paper Closed form (pun intended!) \jour
Contemp. Math.\vol 143\yr 1993\pages 579--607\endref

\ref\key ZPS\by L. L. Zhao, H. Pan and Z. W. Sun\paper Some congruences for the second-order Catalan numbers
\jour Proc. Amer. Math. Soc.\vol 138\yr 2010\pages 37--46\endref

\ref\key Zu\by W. Zudilin\paper Ramanujan-type supercongruences
\jour J. Number Theory\vol 129\yr 2009\pages 1848--1857\endref

\endRefs
\bye